\documentclass[11pt,leqno,openany]{article}
\usepackage{amsmath,amssymb,amsthm,textcomp}
\usepackage{colortbl}
\usepackage[english]{babel}
\usepackage{dsfont}
\newcommand{\bfm}[1]{\boldsymbol{#1}}

\newcounter{ass}

\newtheorem{theorem}{Theorem}[section]
\newtheorem{corollary}[theorem]{Corollary}
\newtheorem{definition}[theorem]{Definition}
\newtheorem{lemma}[theorem]{Lemma}
\newtheorem{proposition}[theorem]{Proposition}

\newtheorem{hypothesis}[ass]{Assumption}
\newtheorem{remark}[theorem]{Remark}

\def \H{ \mathbb H}
\def \M{\mathbb M}
\def \R{ \mathbb R }
\def \reel{ \mathbb R }

\def \one{ {\rm 1}\mkern-4.5mu{\rm I} }

\def \E{ \mathbb  E }

\def \P{ \mathbb P  }

\begin{document}
\title{Critical homogenization of L\'evy process driven SDEs  in random medium}
\author{R\'emi Rhodes${}^{(a)}$ and Ahmadou Bamba  Sow${}^{(b)}$}
\date{}
\maketitle
{\footnotesize \noindent (a)   Universit{\'e} Paris-Dauphine, CEREMADE, Place du Marchal De Lattre de Tassigny, 75775 Paris Cedex 16, France.
Phone: (33)(0)1 44 05 48 51. E-mail: {\tt rhodes@ceremade.dauphine.fr}\\
(b) Universit\'e Gaston Berger, UFR  SAT, LERSTAD,  BP  234, Saint-Louis, S\'en\'egal.\\
Phone : (221) 33 961 23 40. E-mail : {\tt ahbsow@gmail.com}\\}
\begin{abstract}
We are concerned with homogenization of stochastic differential equations (SDE) with stationary coefficients driven by Poisson random measures and Brownian motions in the critical case, that is when the limiting equation admits both a Brownian part as well as a pure jump part.  We state an annealed convergence theorem. This problem is deeply connected with homogenization of integral partial differential equations.
\end{abstract}

\section{Introduction}\label{intro}

If  $B_t$ is a Brownian motion, it is well known that the rescaled process $\epsilon B_{t/\epsilon^2}$ is still a Brownian motion. Starting from this observation, we expect that, under reasonable assumptions on the coefficients, the solution $X$ of the SDE
$$ X_t=x+\int_0^tb(X_r)\,dr+\int_0^t\sigma(X_r)\,dB_r$$ admits a scaling limit, namely that the rescaled process $\epsilon X_{t/\epsilon^2} $ should converge towards a Brownian motion. This problem has been widely studied when the coefficients are periodic or, more recently, when the coefficients are stationary random fields. Quoting all the references is beyond the scope of the paper.

We can make the same observation concerning an $\alpha$-stable L\'evy process $L$: the process $L_t$ and the rescaled process $\epsilon L_{t/\epsilon^\alpha}$ have the same law. This leads to studying scaling limits of SDEs driven by Poisson random measures and, more generally, SDEs driven by both Brownian motions and Poisson random measures (called It\^o-L\'evy type SDEs). However, that issue has been poorly studied so far: see  \cite{Fra} in the case of SDEs with periodic coefficients only driven by Poisson random measures or \cite{rhodes}  for jump processes arising in the context of boundary problems (the reader may also refer to \cite{schwab} for an insight of analytical methods in the context of periodic coefficients).

In \cite{vargas}, the authors investigate the scaling limits of It\^o-L\'evy type SDEs with stationary random coefficients. They prove that there are two possible limiting behaviours, depending on some integrability condition of the compensator of the Poisson random measure. The limiting equation is either a Brownian motion or an $\alpha $-stable L\'evy process. The first situation arises when the jumps of the Poisson measure are small and thus exhibit a diffusive behaviour  (this latter situation was predictable in the light of the wide literature about random walks in random environment). When the Poisson random measure performs sufficiently long jumps, the jump part overscales the Brownian part. This gives rise to the following question: what is the natural framework to make the limiting equation exhibit both a diffusive part and a jump part? That is the issue we investigate in the present paper.

We further stress that our paper is deeply connected to the issue of homogenizing, as $\epsilon\to 0$,  integral partial differential equations (IPDE) with stationary coefficients of the type
$$\partial_t u^\epsilon(t,x)=\partial_x\big(a(\frac{x}{\epsilon})\partial_x u^\epsilon(t,x)\big)+\lim_{\beta\to 0}\int_{\beta<|z|}\big(u^\epsilon(t,x+z)-u^\epsilon(t,x)\big)\frac{c\big(\frac{x}{\epsilon},\frac{z}{\epsilon}\big)}{|z|^{1+\alpha}}dz$$ with suitable boundary conditions. We will address more precisely that connection (and homogenization) in the case of nonlinear problems in a forthcoming paper.


\section{Statements of the problem}\label{statement}

\subsection{Random medium}
We first introduce the notion of random medium (see e.g.
\cite{jikov}) and the necessary background about random media
\begin{definition}\label{medium}
Let $(\Omega ,{\cal G},\mu )$ be a probability space and
$\left\{\tau_{x};x\in \reel\right\}$ a group of measure
preserving transformations acting ergodically on $\Omega $:

1) $\forall A\in {\cal G},\forall x\in \reel$, $\mu (\tau
_{x}A)=\mu (A)$,

2) If for any $x\in \reel$, $\tau _{x}A=A$ then $\mu (A)=0$ or
$1$,

3) For any measurable function ${\bfm g}$ on $(\Omega ,{\cal
G},\mu)$, the function $(x,\omega )\mapsto {\bfm g}(\tau_x \omega)$
is measurable on $(\reel\times\Omega ,{\cal B}(\reel)\otimes
{\cal G})$.
\end{definition}
The expectation with respect to the random medium is denoted by
${\mathbb M}$. We define as usually the spaces $L^p(\Omega ,{\cal G},\mu )$ for $p\in[1,+\infty]$, or $L^p(\Omega)$ for short. The corresponding norm are denoted by $|\, \cdot \,|_p$. The inner product in $L^2(\Omega)$ is denoted by $(\,
\cdot\,,\, \cdot\,)_2$. The operators on $L^2(\Omega )$ defined by
$T_{x}{\bfm g}(\omega )={\bfm g}(\tau_{x}\omega )$ form a strongly
continuous group of unitary maps in $L^2(\Omega )$. Each function
${\bfm g}$ in $ L^2(\Omega )$ defines in this way a stationary
ergodic random field on $\reel $. The group possesses a
generator $D$, defined by
\begin{equation} D{\bfm g}=\lim_{\R\ni h\to
0}h^{-1}(T_{h}{\bfm g}-{\bfm g}) \text{ if the limit exists in
the }L^2(\Omega)\text{-sense},
\end{equation}
which is closed and densely defined. We distinguish the differential operator in random medium $D$ from the usual derivative $\partial_x f $ of a function $f$ defined on $\R$.

\noindent {\bf Notations.} Recursively, we define the operators ($k\geq 1$) $D^k=D(D^{k-1})$ with domain $H^k(\Omega)=\{{\bfm f}\in H^{k-1}(\Omega); D^{k-1}{\bfm f}\in {\rm Dom}({\bfm D})=H^1(\Omega)\}$. We also define $H^\infty(\Omega)=\bigcap_{k=1}^\infty H^k(\Omega)$.

We denote with ${\cal C}$ the dense subspace of $L^2(\Omega )$ defined by
\begin{equation*}
{\cal C}={\rm Span}\left\{{\bfm g} \star \varphi ;{\bfm g}\in
L^\infty(\Omega ),\varphi \in C^\infty _c(\reel  )\right\}
\end{equation*}
with ${\bfm g} \star \varphi(\omega )=\int_{\reel }{\bfm
g}(\tau_{x}\omega )\varphi (x)\,dx$. We point out that ${\cal
C}\subset {\rm Dom}(D)$, and $D({\bfm
g} \star \varphi)=-{\bfm g} \star
\partial  \varphi/\partial x$. This last quantity is also equal to $D{\bfm g}
\star \varphi $ if ${\bfm g}\in {\rm Dom}(D)$. $C(\Omega)$ is defined as the closure of $\mathcal{C} $ in $L^\infty(\Omega)$ with respect to the norm $|\cdot|_\infty$, whereas $C^\infty(\Omega)$ stands for the subspace of $H^\infty(\Omega)$, whose elements satisfy: ${\bfm f}\in C^\infty(\Omega)\Leftrightarrow \forall k \geq  0,\,|D^k{\bfm f} |_\infty<+\infty$. We point out that, whenever a function ${\bfm f}\in H^\infty(\Omega)$, $ \mu$
a.s. the mapping $f_\omega:x\in\R\mapsto {\bfm f}(\tau_x\omega)$  is infinitely differentiable and $\partial_x f_\omega(x)=D{\bfm f}(\tau_x\omega)$.
\subsection{Structure of the coefficients}
We consider a  L\'evy measure $ \nu$, that is a  $\sigma$-finite measure  $\nu$ on $\R$ of the type
\begin{equation}\label{nu}
\nu(dz)=\frac{1}{|z|^{1+\alpha}}dz
\end{equation}
for some $\alpha\in]0,2[$.
We introduce the coefficients $ {\bfm V},{\bfm \sigma}\in L^\infty(\Omega)$ and $ {\bfm \gamma}:\Omega\times\R\to \R$ satisfying the following conditions.

For each fixed $\omega\in \Omega$, by defining  the mapping $\gamma_\omega:z\mapsto {\bfm \gamma}(\omega,z)$, we can consider the measure $\nu\circ  \gamma_\omega^{-1}:A\subset \R\mapsto \nu(\gamma_\omega^{-1} (A))=\nu\big(\{z\in\R; {\bfm \gamma}(\omega,z)\in A\}\big)$.

\begin{hypothesis}\label{symmetry}{\bf Symmetry of the kernel.}
We assume that the measure $\nu\circ  \gamma_\omega^{-1}$ can be rewritten as
$$\nu\circ  \gamma_\omega^{-1}(dz)=e^{2{\bfm V}(\omega)}\frac{{\bfm c}(\omega,z)}{|z|^{1+\alpha}}\,dz$$ for some  measurable nonnegative symmetric kernel ${\bfm c}$ defined on $\Omega\times \R$. The symmetry of ${\bfm c}$ means
$$\mu \text{ a.s.},\,\,dz \,\text{ a.s.}, \quad {\bfm c}(\tau_z\omega,-z)={\bfm c}(\omega,z).\qed$$
 \end{hypothesis}

 \begin{hypothesis}\label{limitc}{\bf Limiting kernel.}
 We assume that there exists a function ${\bfm \theta}\in L^1(\Omega)$ such that
 $$\lim_{|z|\to\infty}\M\big[|{\bfm c}(\omega,z)-{\bfm \theta}(\omega)|\big]=0. $$
 \end{hypothesis}

 \begin{hypothesis}\label{ellipticity}{\bf Ellipticity.} We set ${\bfm a}={\bfm \sigma}^2$. There is a constant $M_{\ref{ellipticity}}>0$ such that
$$\forall (\omega,z)\in \Omega\times \R,\quad M_{\ref{ellipticity}}^{-1}\leq {\bfm a}(\omega)\leq M_{\ref{ellipticity}}\quad \text{ and }\quad M_{\ref{ellipticity}}^{-1}\leq {\bfm c}(\omega,z)\leq M_{\ref{ellipticity}}.\qed$$
\end{hypothesis}

\begin{hypothesis}\label{regul}{\bf Regularity.} We assume the coefficients satisfy the following assumptions:

1) The coefficients ${\bfm V}$, ${\bfm \sigma}$ belong to $C^\infty(\Omega)$. In particular, we can define $${\bfm b}=\frac{1}{2}D{\bfm a}-{\bfm a}D{\bfm V} =\frac{e^{2{\bfm V} }}{2}D\big(e^{-2{\bfm V}}{\bfm a} \big)\in C^\infty(\Omega),$$

2) For  $dz$-almost every $z\in\R$, the mapping $\omega\mapsto {\bfm c}(\omega,z)$ belongs to $C^\infty(\Omega)$ and, for each fixed $k\geq 1$, there exists a constant $C_k$ such that $|D^k{\bfm c}(\cdot,z)|_\infty\leq C_k$, $dz$ a.s.

3) $\mu$ a.s., for $dz$ almost every $|z|>1$, the mapping $x\in\R\mapsto {\bfm \gamma}(\tau_x\omega,z)$ is continuous and $\mu$ a.s., we can find a constant $C>0$ such that $\forall x,y\in\R $, $$\int_{|z|\leq 1}|{\bfm \gamma}(\tau_y\omega,z)-{\bfm \gamma}(\tau_x\omega,z)|^2\nu(dz)\leq C |y-x|^2.$$

4) For every  $\omega\in  \Omega$, the limit $${\bfm
e}(\omega)=\lim_{\beta\to 0}\int_{\beta\leq |{\bfm
\gamma}(\omega,z)|}{\bfm \gamma}(\omega,z)\one_{|z|\leq1}\nu(dz)$$
exists in the $ L^2(\Omega)$ sense and defines a bounded
Lipschitzian function, that is (for some constant
$M_{\ref{regul}}\geq 0$), $|{\bfm e}|_\infty\leq M_{\ref{regul}}$
and $ \mu$  a.s., $\forall x,y\in\R$, $ |{\bfm
e}(\tau_y\omega)-{\bfm e}(\tau_x\omega)|\leq M_{\ref{regul}}|x-y|$.

5) Furthermore, there is a positive constant $S$ such that $\sup_{|z|\leq 1}|{\bfm \gamma}(\cdot,z)|_{\infty}\leq S$
\qed
 \end{hypothesis}


Even if it means adding to ${\bfm V}$ a renormalization constant (this does not change the drift ${\bfm b}$ and the jump coefficients ${\bfm \gamma} $ and $\nu$), we may assume that 
$$\M[e^{-2{\bfm V}}]=1 $$
and consider the probability measure $d\pi=e^{-2{\bfm V}}\,d\mu$ on $(\Omega,\mathcal{G})$. We denote by $\M_\pi$ the expectation w.r.t. this probability measure.

\subsection{Jump-diffusion processes in random medium}
 We suppose that we are
given a complete probability space  $(\Omega',{\cal F},\P)$ with a
right-continuous increasing family of complete sub $\sigma$-fields
$({\cal F}_t)_t$ of ${\cal F}$, a ${\cal F}_t$-adapted Brownian
motion $ \{B_t;t\geq 0\}$ and ${\cal F}_t$-adapted Poisson random
measure $N(dt,dz)$ with intensity $\nu$.
$\widetilde{N}(dt,dz)=N(dt,dz)-\nu(dz)dt$ denotes the corresponding
compensated random measure and $\hat{N}^\epsilon(dt,dz)$ the
truncated compensated random measure $N(dt,dz)-\one_{|z|\leq
\epsilon}\nu(dz)dt, \;\epsilon>0$. We further assume that the
Brownian motion, the Poisson random measure and the random medium are
independent.

For each fixed $ \omega\in\Omega$ and $\epsilon>0$, Assumptions \ref{regul}.3 and \ref{regul}.4 are enough to ensure  existence and pathwise uniqueness of a ${\cal F}_t$-adapted process $X^\epsilon$ (see \cite[Ch.6, Sect.2]{applebaum}) solution to the following SDE
\begin{equation}\label{SDE}
\begin{split}
X^\epsilon_t=x+\int_0^t (\frac{1}{\epsilon}{\bfm b}+
\epsilon^{1-\alpha}{\bfm
e})(\tau_{\overline{X}^\epsilon_{r-}}\omega)\,dr+\int_0^t
\int_{\R}\epsilon{\bfm
\gamma}(\tau_{\overline{X}^\epsilon_{r-}}\omega,\frac{z}{\epsilon})\,\hat{N}^\epsilon(dr,dz)+\int_0^t{\bfm
\sigma}(\tau_{\overline{X}^\epsilon_{r-}}\omega)\,dB_r,
\end{split}
\end{equation}
where we have set $\overline{X}^\epsilon_t= X^\epsilon_t/\epsilon$.

\begin{remark}
The reader may find the formulation of SDE \eqref{SDE} a bit disturbing. Actually, this is just a correct formulation for the SDE we have in mind, namely
$$X^\epsilon_t=x+\int_0^t \frac{1}{\epsilon}{\bfm b}(\tau_{\overline{X}^\epsilon_{r-}}\omega)\,dr+\int_0^t \int_{\R}\epsilon{\bfm \gamma}(\tau_{\overline{X}^\epsilon_{r-}}\omega,\frac{z}{\epsilon})\,N(dr,dz)+\int_0^t{\bfm \sigma}(\tau_{\overline{X}^\epsilon_{r-}}\omega)\,dB_r. $$
Due to integrability issues, the above formal equation admits the correct formulation \eqref{SDE}.
\end{remark}
\subsection{Main result}
We denote with  $D(\R_+;\R)$  the space of right-continuous $\R$-valued functions with left limits, endowed with the Skorohod topology, cf \cite{ethier}. We fix $x\in \R$ and we claim
\begin{theorem}\label{mainresult}
In $\mu$ probability, the process $X^\epsilon$, starting from $x\in \R$, converges in law in the Skorohod topology towards a L\'evy process $L$ with characteristic function $\E[e^{iu L_t}]=e^{t \varphi(u)}$, where the L\'evy exponent $\varphi$ is given by
$$\varphi(u)=-\frac{1}{2}Au^2+\M[{\bfm \theta}]\int_{\R}(e^{iuz}-1-i u z\one_{|z|\leq 1})\nu(dz) $$ for some constant coefficient $A$.
\end{theorem}

\begin{remark}
Actually, by looking closely in the proofs in Section \ref{sec:correctors}, we could prove that $A$ exactly matches the homogenized coefficient when the SDE \eqref{SDE} possesses no jump part. In particular, we could prove the variational formula
\begin{equation}\label{variational}
A=\inf_{{\bfm \varphi}\in \mathcal{C}} \M_\pi\big[{\bfm a}(1+D{\bfm \varphi})^2\big],
\end{equation}
from which lower and upper bounds for $A$ can be obtained. It is then plain to see that $A$ is nondegenerate (because ${\bfm a}$ is, see \cite{olla} for the derivation of the variational formula).
\end{remark}

\section{Dirichlet forms in random medium}\label{sec:dirichlet}

We can equip the space $L^2(\Omega)$ with the inner product $({\bfm \varphi},{\bfm \psi})_\pi=\M[{\bfm \varphi}{\bfm \psi}e^{-2{\bfm V}}]$, and denote by $|\cdot|_\pi$ the associated norm. Since ${\bfm V}$ is bounded, both inner products $(\cdot,\cdot)_2$ and $ (\cdot,\cdot)_\pi$ are equivalent on $L^2(\Omega)$.


Let us define on ${\cal C}\times {\cal C} $ the following bilinear forms (with $\lambda>0$)
\begin{equation}\label{dirichletf}
\begin{split}
B^d({\bfm \varphi},{\bfm \psi})&=\frac{1}{2}({\bfm a}D{\bfm \varphi},D{\bfm \psi})_\pi,\\ B^{j}({\bfm \varphi},{\bfm \psi})&=\frac{1}{2}\M\int_{\R}(T_z{\bfm \varphi}-{\bfm \varphi})(T_z{\bfm \psi}-{\bfm \psi}){\bfm c}(\cdot, z)\nu(dz),\\
B_\lambda^\epsilon({\bfm \varphi},{\bfm \psi})&=\lambda({\bfm \varphi},{\bfm \psi})_\pi+B^d({\bfm \varphi},{\bfm \psi})+\epsilon^{2-\alpha}B^{j}({\bfm \varphi},{\bfm \psi})
\end{split}
\end{equation}
We can thus consider on ${\cal C}\times {\cal C}  $ the inner product $B_\lambda^\epsilon $ and the closure $\H$ of ${\cal C} $ w.r.t. the associated norm (note that the definition of $\H$ does not depend on $\lambda,\epsilon>0$ since the corresponding norms are equivalent from Assumption \ref{ellipticity}).

Fix $\epsilon>0$. In what follows, we use the same strategy as in \cite[Sect. 3]{vargas} to which  the reader is
referred for further details (as well as the references therein). For any $\lambda>0$, $B_\lambda^\epsilon$ extends to $\H\times \H$.
This extension, still denoted $B_\lambda^\epsilon$,  defines a resolvent operator $G^\epsilon_\lambda: L^2(\Omega)\to \H$, which is
one-to-one and continuous. It thus defines an unbounded operator ${\bfm L}^\epsilon=\lambda-(G^\epsilon_\lambda)^{-1}$ on $L^2(\Omega)$
with domain ${\rm Dom}({\bfm L}^\epsilon)=G_\lambda^\epsilon(L^2(\Omega))$. This definition does not depend on $\lambda>0$.
The unbounded operator
${\bfm L}^\epsilon$ is closed, densely defined and seld-adjoint. We further stress that the weak form of the resolvent
equation $\lambda G^\epsilon_\lambda{\bfm f}-{\bfm L}^\epsilon G^\epsilon_\lambda{\bfm f}={\bfm f} $ reads: $\forall {\bfm \psi}\in\H$
\begin{align}
\lambda(G^\epsilon_\lambda{\bfm f},{\bfm \psi})_\pi+&\frac{1}{2}({\bfm a}DG^\epsilon_\lambda{\bfm f},D{\bfm \psi})_\pi+\frac{\epsilon^{2-\alpha}}{2}\M\int_{\R}(T_zG^\epsilon_\lambda{\bfm f}-G^\epsilon_\lambda{\bfm f})(T_z{\bfm \psi}-{\bfm \psi}){\bfm c}(\cdot, z)\nu(dz)\nonumber\\\label{eqbase}&= ({\bfm f},{\bfm \psi})_\pi.
\end{align}

For sufficiently smooth functions, ${\bfm L}^\epsilon$ can be easily identified (the proof does not differ from \cite[Lemma 3.1]{vargas}): if ${\bfm \varphi}\in H^\infty(\Omega)$,
then ${\bfm \varphi}\in{\rm Dom}({\bfm L}^\epsilon)$ and
\begin{align}\label{generator}
 {\bfm L}^\epsilon{\bfm \varphi}&=\epsilon^{2-\alpha}{\bfm e}D{\bfm \varphi}+\epsilon^{2-\alpha}\int_{\R}\big({\bfm \varphi}(\tau_{{\bfm \gamma}(\omega,z)}\omega)- {\bfm \varphi}(\omega)-{\bfm \gamma}(\omega,z)\one_{\{|z|\leq 1\}}D{\bfm \varphi}(\omega)\big)\,\nu(dz)\\
 &+\frac{1}{2}{\bfm a}D^2{\bfm \varphi}+{\bfm b}D{\bfm \varphi}\nonumber.
\end{align}

Following the proof in \cite{vargas}, we can prove:
\begin{proposition}\label{core}
1) For each $\lambda>0$, the resolvent operator $G^\epsilon_\lambda$ maps $L^2$ into $H^2(\Omega)$, and $H^m(\Omega)$ into $H^{m+2}(\Omega)$ for any $m\geq 1$. In particular ${\rm Dom}(({\bfm L}^\epsilon)^m)=H^{2m}(\Omega) $.

2) The self-adjoint operator ${\bfm L}^\epsilon$ generates a strongly continuous contraction semi-group $(P^\epsilon_t)_t$ of self-adjoint operators. Moreover, we have
 \begin{align}\label{reg:semigroup}
 {\bfm f}\in L^2(\Omega)&\Rightarrow t\mapsto P^\epsilon_t{\bfm f} \in C([0;+\infty[;L^2(\Omega))\cap C^\infty(]0;+\infty[;H^\infty(\Omega)),\\
{\bfm f}\in H^\infty(\Omega)&\Rightarrow t\mapsto P^\epsilon_t{\bfm f} \in C^\infty([0;+\infty[;H^\infty(\Omega))
\end{align}
where, given an interval $I\subset\R$, $C(I;L^2(\Omega))$ (resp. $C^\infty(I;H^\infty(\Omega))$) stands for the space of continuous functions from $I$ to $L^2(\Omega)$ (resp. infinitely differentiable functions from $I$ to $ H^\infty(\Omega)$).

3) The semi-group $(P^\epsilon_t)_t$ is sub-Markovian. Put in other
words, for any ${\bfm f}\in L^2(\Omega)$ such that $0\leq {\bfm
f}\leq 1$ $\mu$ a.s., we have $0\leq P^\epsilon_t{\bfm f}\leq 1$
$\mu$ a.s. for any $t>0$. In particular,
$P^\epsilon_t:L^\infty(\Omega)\to L^\infty(\Omega)$ and
$G^\epsilon_\lambda: L^\infty(\Omega)\to L^\infty(\Omega)$ are
continuous with respective norms $1$ and $1/\lambda$.
\end{proposition}

Similarly, we can consider on $\mathcal{C}\times \mathcal{C}$ the bilinear form
$$D_\lambda({\bfm \varphi},{\bfm \psi})=\lambda({\bfm \varphi},{\bfm \psi})_\pi+B^d({\bfm \varphi},{\bfm \psi}).$$
The form $D_\lambda$ on ${\cal C}\times {\cal C}  $ defines an inner product  and we can define the closure $\H^d $ of ${\cal C}\times {\cal C}  $ with respect to $ D_\lambda$. Once again, the definition of $\H^d$ does not depend on $\lambda>0$ since the corresponding norms are equivalent.
For any $\lambda>0$, $D_\lambda$ continuously extends to $\H^d\times \H^d$. This extension is still denoted $D_\lambda$ and  defines a resolvent operator $G^d_\lambda: L^2(\Omega)\to \H^d$, which is one-to-one and continuous. It thus defines an unbounded operator ${\bfm L}^d=\lambda-(G^d_\lambda)^{-1}$ on $L^2(\Omega)$ with domain ${\rm Dom}({\bfm L}^d)=G_\lambda^d(L^2(\Omega))$. This definition does not depend on $\lambda>0$. Moreover, ${\bfm L}^d$ is self-adjoint and it is plain to see that ${\bfm L}^d$ is given on $\mathcal{C}$ by
$${\bfm L}^d{\bfm \varphi}=\frac{e^{2{\bfm V}}}{2}D\big({\bfm a}e^{-2{\bfm V}}D{\bfm \varphi}\big)=\frac{1}{2}{\bfm a}D^2{\bfm \varphi}+{\bfm b}D{\bfm \varphi}.$$

\begin{lemma}\label{lem:s}
If the function ${\bfm \varphi}$ belongs to  ${\rm Dom}({\bfm L}^d)$, then  ${\bfm \varphi}$ also belongs to $\H$.
\end{lemma}

\noindent {\it Proof of Lemma \ref{lem:s}.} It is plain to see that the lemma results from the following inequality
$$\forall \bfm \varphi\in \mathcal{C},\quad B^j( \bfm \varphi, \bfm \varphi)\leq C D_\lambda( \bfm \varphi, \bfm \varphi)$$ for some positive constant $C$ that may depend on $\lambda$. Since this result is quite classical, details are left to the reader.\qed

\section{Invariant measure}\label{measure}
In what follows, $X^\epsilon$ denotes the solution of \eqref{SDE}
starting from $0$.

\begin{proposition}\label{prop:invar}
For each function $ {\bfm f}\in C(\Omega)$, we have $$\M_\pi[\E[{\bfm f}(\tau_{\overline{X}^\epsilon_t}\omega)]]=\M_\pi[\E[{\bfm f}(\tau_{\overline{X}^\epsilon_{t-}}\omega)]]= \M_\pi[{\bfm f}].$$
\end{proposition}

\noindent {\it Proof.} Given ${\bfm \varphi}\in \mathcal{C}\subset H^\infty(\Omega)$ and $t>0$, the mapping $(s,\omega)\mapsto P^\epsilon_{\epsilon^{-2}(t-s)}{\bfm \varphi}$ belongs to $ C^\infty([0,t]; H^\infty( \Omega))$ and is bounded (cf Prop \ref{core}). We can thus apply the It\^o formula between $ 0$ and $t$, which reads  $\mu$ a.s. (use $\partial_tP_t^\epsilon{\bfm \varphi}={\bfm L}^\epsilon P_t^\epsilon{\bfm \varphi}$):
\begin{align}\label{tim}
{\bfm
\varphi}&(\tau_{\overline{X}^\epsilon_{t}}\omega)=P^\epsilon_{\epsilon^{-2}t}{\bfm
\varphi}(\omega)+\epsilon^{-1}\int_0^tDP^\epsilon_{\epsilon^{-2}(t-r)}{\bfm
\varphi}{\bfm
\sigma}(\tau_{\overline{X}^\epsilon_{r-}}\omega)\,dB_r\\&+\int_0^t\big(P^\epsilon_{\epsilon^{-2}(t-r)}{\bfm
\varphi}(\tau_{\overline{X}^\epsilon_{r-}+{\bfm
\gamma}(\tau_{\overline{X}^\epsilon_{r-}}\omega,z)}\omega)-P^\epsilon_{\epsilon^{-2}(t-r)}{\bfm
\varphi}(\tau_{\overline{X}^\epsilon_{r-}}\omega)\big)\,\widetilde{N}(dr,dz)\nonumber
\end{align}
We remind the reader that $\mu$ a.s., $\P( X^\epsilon\, \text{ is c\`ad-l\`ag on }[0,+\infty[)=1$. Hence, for any $t>0$, we have  $\P(\sup_{0\leq s \leq t}|X^\epsilon_s|<+\infty)=1$. We deduce that the sequence of stopping times $ S_n=\inf\{s\geq 0; |X^\epsilon_s|>n\}$ satisfies: $\mu$ a.s., $\P$ a.s. $S_n\to +\infty$ as $n\to \infty$. By replacing $t$ by $t\wedge S_n$ (i.e. $\min(t,S_n) $) in \eqref{tim} and by taking the expectation, the martingale terms vanish and we get
$$\E[{\bfm \varphi}(\tau_{\overline{X}^\epsilon_{t\wedge S_n}}\omega)] =\E[P^\epsilon_{\epsilon^{-2}t\wedge S_n}{\bfm \varphi}(\omega)].$$
Using the boundedness of ${\bfm \varphi}$, $P^\epsilon_t{\bfm \varphi}$ and the continuity of the mappings $x\mapsto {\bfm \varphi}(\tau_x\omega)$, $t\mapsto  P^\epsilon_t{\bfm \varphi}$, we can pass to the limit as $n\to \infty$ in the above equality  to prove $P^\epsilon_{\epsilon^{-2}t} {\bfm \varphi}(\omega)=\E[{\bfm \varphi}(\tau_{\overline{X}^\epsilon_{t}}\omega)]$ for $  {\bfm \varphi}\in \mathcal{C}$.

In case ${\bfm f}\in C(\Omega)$, we can find a sequence $({\bfm \varphi}_n)_n \in \mathcal{C}$ converging towards ${\bfm f}$ in $L^\infty(\Omega)$-norm (for instance $({\bfm f}\star \rho_n)_n$ for some regularizing sequence $(\rho_n)_n\subset C^\infty_c(\R)$). By passing to the limit in the relation $P^\epsilon_{\epsilon^{-2}t} {\bfm \varphi}_n(\omega)=\E[{\bfm \varphi}_n(\tau_{\overline{X}^\epsilon_{t}}\omega)]$, we get the relation $P^\epsilon_{\epsilon^{-2}t} {\bfm f}(\omega)=\E[{\bfm f}(\tau_{\overline{X}^\epsilon_{t}}\omega)]$ for each ${\bfm f}\in C(\Omega)$. Finally, we complete the proof by noticing that $\M_\pi[P^\epsilon_t{\bfm f}]=\M_\pi[{\bfm f}]$ by construction of $(P^\epsilon_t)_t$ and that $\E[{\bfm f}(\tau_{\overline{X}^\epsilon_{t}}\omega)]=\E[{\bfm f}(\tau_{\overline{X}^\epsilon_{t-}}\omega)]$ since $\forall t$, $\P\big(X^\epsilon_t=X^\epsilon_{t-}\big)=1$.\qed

\begin{remark}
As a direct consequence, for each $\epsilon>0$, the mapping $${\bfm f}\in C(\Omega)\mapsto \big(t\mapsto \int_0^t{\bfm f}(\tau_{\overline{X}^\epsilon_{r}}\omega)\,dr\in L^1(\Omega;C([0,T]))\big)$$ continuously (and uniquely) extends to $L^1(\Omega)$. The extension is still denoted by $ \int_0^t{\bfm f}(\tau_{\overline{X}^\epsilon_{r}}\omega)\,dr$ for ${\bfm f}\in L^1(\Omega)$.
\end{remark}
\section{Ergodic problems}

The main purpose of this section is to establish the following  results:
\begin{theorem}{\bf Ergodic theorem I.}\label{ergth}
For any ${\bfm f}\in L^1(\Omega)$, the following convergence holds
$$\lim_{\epsilon\to 0}\M_\pi\E\Big[\sup_{0\leq t \leq T}\big|\int_0^t{\bfm f}(\tau_{\overline{X}^\epsilon_{r-}}\omega)\,dr-t\M_\pi[{\bfm f}]\big|\Big]=0. $$
\end{theorem}

\begin{corollary}\label{corerg2}
Given a family $ ({\bfm f}_\epsilon)_\epsilon$ converging towards $ {\bfm f}\in L^1(\Omega)$, we have
$$\lim_{\epsilon\to 0}\M_\pi\E\Big[\sup_{0\leq t \leq T}\big|\int_0^t{\bfm f}_\epsilon(\tau_{\overline{X}^\epsilon_{r-}}\omega)\,dr-t\M_\pi[{\bfm f}]\big|\Big]=0. $$
\end{corollary}

\begin{corollary}\label{corerg}
Given a continuous  function $g$ satisfying $|g(z)|\leq \min(1,z^2)$ for all $z\in\R$, the following convergence holds
$$\lim_{\epsilon\to 0}\M_\pi\E\Big[\big|\int_0^t\int_\R g\big(\epsilon{\bfm \gamma}(\tau_{\overline{X}^\epsilon_{r-}}\omega,\frac{z}{\epsilon})\big) \nu(dz)dr-t\M[{\bfm \theta}]\int_\R g(z) \nu(dz)\big|\Big]=0.$$
\end{corollary}

\noindent {\it Proof of Corollary \ref{corerg2}.} Clearly, the result follows from Theorem \ref{ergth} applied to ${\bfm f} $ and from the inequality
$$\M_\pi\E\Big[\sup_{0\leq t \leq T}\big|\int_0^t({\bfm f}_\epsilon-{\bfm f})(\tau_{\overline{X}^\epsilon_{r-}}\omega)\,dr\big|\Big]\leq T\M_\pi[|{\bfm f}_\epsilon-{\bfm f}|].$$
\qed

\noindent {\it Proof of Corollary \ref{corerg}.} Consider such a function $g$. We have
\begin{align*}
\int_0^t\int_\R g\big(\epsilon{\bfm \gamma}(\tau_{\overline{X}^\epsilon_{r-}}\omega,\frac{z}{\epsilon})\big)\nu(dz)dr=&\epsilon^{-\alpha}\int_0^t\int_\R g\big(\epsilon{\bfm \gamma}(\tau_{\overline{X}^\epsilon_{r-}}\omega,z)\big)\nu(dz)dr\\
=&\epsilon^{-\alpha}\int_0^t\int_\R g(\epsilon z){\bfm c}(\tau_{\overline{X}^\epsilon_{r-}}\omega,z)e^{2{\bfm V}(\tau_{\overline{X}^\epsilon_{r-}}\omega)}\nu(dz)dr\\
=&\int_0^t\int_\R g(z){\bfm
c}(\tau_{\overline{X}^\epsilon_{r-}}\omega,\frac{z}{\epsilon})e^{2{\bfm
V}(\tau_{\overline{X}^\epsilon_{r-}}\omega)}\nu(dz)dr.
\end{align*}
We define ${\bfm G}^\epsilon(\omega)=\int_\R g(z){\bfm
c}(\omega,\frac{z}{\epsilon})e^{2{\bfm
V}(\omega)}\nu(dz)$. From Assumption \ref{limitc}, the family $({\bfm G}^\epsilon)_\epsilon$ converges towards ${\bfm G}(\omega)={\bfm
\theta}(\omega)e^{2{\bfm
V}(\omega)}\int_\R g(z)\nu(dz)$ in $L^1(\Omega)$. 
The proof can be completed with Corollary \ref{corerg2}.\qed

The proof of Theorem \ref{ergth} is based on several auxiliary results listed (and proved) below
\begin{proposition}\label{prop:erg}
Given ${\bfm f}\in L^2(\Omega)$ and a family
$(\lambda(\epsilon))_{\epsilon>0}\subset ]0;+\infty[$ such that
$$\lim_{\epsilon\to 0}\lambda(\epsilon)=0,$$ we define ${\bfm
u}^\epsilon$ (for any $\epsilon>0$) as the solution of the resolvent
equation
$$\lambda(\epsilon) {\bfm u}^\epsilon-{\bfm L}^\epsilon{\bfm u}^\epsilon = {\bfm f}.$$
Then we have
$$\lim_{\epsilon\to 0}|\lambda(\epsilon) {\bfm u}^\epsilon -\M_\pi[{\bfm
f}]|_2=0$$ and the estimates
\begin{equation}\label{est:erg}
\lambda(\epsilon)|D{\bfm u}^\epsilon|_\pi^2+\frac{\lambda(\epsilon)\epsilon^{2-\alpha}}{2}\M\int_\R
(T_z {\bfm u}^\epsilon - {\bfm u}^\epsilon)^2 {\bfm c}(\omega,
z)\nu(dz)\leq |\bfm f|^2_\pi.
\end{equation}
\end{proposition}

\noindent {\it Proof.} From  the  resolvent  equation
\eqref{eqbase} (where we choose $ {\bfm \psi}={\bfm u}^\epsilon$),   we have
\begin{align}
\lambda(\epsilon)  |\bfm  u^\epsilon|_\pi^2  + \frac{1}{2}({\bfm a}
D {\bfm u^\epsilon}, D {\bfm  u^\epsilon})_\pi &+
\frac{\epsilon^{2-\alpha}}{2} \M\int_\R (T_z {\bfm u}^\epsilon -
{\bfm u}^\epsilon)^2 {\bfm c}(\omega, z)\nu(dz) =
({\bfm  f},{\bfm u}^\epsilon)_\pi. \label{multiplie-lambda} 
\end{align}
We use the  Cauchy-Schwarz   inequality   in  the   right-hand   side to obtain:
$$
({\bfm  f},{\bfm u}^\epsilon)_\pi   \le  \frac{|{\bfm  f}|_\pi^2}{2
\lambda(\epsilon)}     + \frac{\lambda(\epsilon)  |\bfm
u^\epsilon|_\pi^2}{2}.
$$
By plugging   this  inequality  into  \eqref{multiplie-lambda} and by multiplying by $\lambda(\epsilon)$,  we
obtain
\begin{align}
\frac{\lambda^2(\epsilon)  |\bfm  u^\epsilon|_\pi^2}{2}  +
\frac{\lambda(\epsilon)}{2}({\bfm a} D {\bfm u^\epsilon}, D {\bfm
u^\epsilon})_\pi &+ \frac{\lambda(\epsilon)}{2}\epsilon^{2-\alpha}
\M\int_\R (T_z {\bfm u}^\epsilon - {\bfm u}^\epsilon)^2 {\bfm
c}(\omega, z)\nu(dz) \; \le \frac{|{\bfm  f}|_\pi^2}{2}.
\end{align}
Estimate \eqref{est:erg} then results from Assumption
\ref{ellipticity}. Moreover, the family $(\lambda(\epsilon)  \bfm
u_\epsilon)_\epsilon$ is  bounded in $L^2(\Omega)$ and,  along a
subsequence, we can find $\overline{\bfm f}\in L^2(\Omega)$
such that
$$
\lambda(\epsilon) \bfm  u^\epsilon \rightarrow
\overline{{\bfm f}}\quad \text{weakly in }L^2(\Omega).
$$
From the   resolvent equation \eqref{eqbase}, we  have for any
$\bfm \varphi \in \H$,
\begin{align}
\lambda(\epsilon)^{2}(\bfm  u^\epsilon, \bfm \varphi)_\pi   +
\frac{\lambda(\epsilon)}{2}({\bfm a} D {\bfm u^\epsilon}, D
\bfm\varphi)_\pi &+
\frac{\lambda(\epsilon)}{2}\epsilon^{2-\alpha}\M\int_\R (T_z {\bfm
u}^\epsilon - {\bfm u}^\epsilon) (T_z \bfm\varphi - \bfm\varphi)
{\bfm c}(\omega, z)\nu(dz) \nonumber\\
&= ({\bfm f},\bfm\varphi)_\pi \lambda(\epsilon).\label{eq:rescor-1}
\end{align}
Thanks to  lemma \ref{lem:s},  \eqref{eq:rescor-1} also  holds   for any
$ \bfm \varphi \in {\rm Dom}({\bfm L}^d)$.

We now investigate the limit of each term in \eqref{eq:rescor-1} as
$ \epsilon\to   0$ when the function ${\bfm \varphi} $ is assumed to belong to $ {\rm Dom}({\bfm L}^d)$. We use the relation (valid for $\bfm\varphi \in
{\rm Dom}({\bfm L}^d)$):  $ \frac{1}{2}({\bfm a} D {\bfm
u^\epsilon}, D \bfm\varphi)_\pi =   - ({\bfm u^\epsilon}, {\bfm L}^d
\bfm \varphi)_\pi $ and we deduce
$$
\lambda(\epsilon)({\bfm a} D {\bfm u^\epsilon}, D \bfm\varphi)_\pi
 =   - (\lambda(\epsilon){\bfm u^\epsilon}, {\bfm L}^d\varphi)_\pi
\rightarrow   - (\overline{{\bfm  f}}, {\bfm L}^d
\bfm \varphi)_\pi\quad \text{ as }\epsilon\to 0. $$ On the other hand, from \eqref{est:erg}, we
have $$ \lim_{\epsilon\to 0} \Big[ \lambda(\epsilon)^{2}(\bfm
u^\epsilon, \bfm \varphi)_\pi  +
\frac{\lambda(\epsilon)}{2}\epsilon^{2-\alpha}\M\int_\R (T_z {\bfm
u}^\epsilon - {\bfm u}^\epsilon) (T_z \bfm\varphi - \bfm\varphi)
{\bfm c}(\omega, z)\nu(dz) \Big] \leq  \lim_{\epsilon\to 0} \lambda(\epsilon)|{\bfm f}|_\pi^2/2=  0.$$ So we are in position to
pass to the limit as  $ \epsilon\to   0$  in \eqref{eq:rescor-1} and
we obtain $(\overline{{\bfm f}}, {\bfm L}^d \bfm\varphi)_\pi = 0$
for any $\bfm \varphi \in {\rm Dom}({\bfm L}^d)$.  Since ${\bfm
L}^d$ is self-adjoint, we deduce $\overline{{\bfm  f}}\in  {\rm
Dom}({\bfm L}^d)$ and ${\bfm L}^d \overline{{\bfm f}}=0$. In
particular, $({\bfm L}^d\overline{{\bfm  f}}, \overline{{\bfm  f}})
= - \frac{1}{2}({\bfm a}  D \overline{{\bfm  f}}, D\overline{{\bfm
f}})=  0$. As a consequence we deduce that $\overline{{\bfm f}}$ is
constant $\mu$ almost surely.

We now determine the constant $\overline{{\bfm  f}}$.
Plugging   the   function  $ \bfm \psi=  \bfm 1$  into the relation
(\ref{eqbase}) yields for  every $\epsilon>0$,
$$\M_\pi[\lambda(\epsilon){\bfm u^\epsilon}]  =  \M_\pi[{\bfm f}].$$
It just remains to let  $\epsilon$ go to $ 0$ and to use  the  weak
convergence to obtain $ \M_\pi[\overline{{\bfm  f}}]  = \M_\pi[{\bfm f}]$. As
a consequence $ \M_\pi[{\bfm f}]   = \overline{{\bfm  f}}
$. This establishes uniqueness of the weak limit of the family $(\lambda(\epsilon){\bfm u^\epsilon})_\epsilon$ so that the whole family is weakly converging (not along a subsequence).

Finally we establish the strong convergence of the family
$(\lambda(\epsilon){\bfm u^\epsilon})_\epsilon$ towards $
\overline{{\bfm  f}}$ as $\epsilon\to 0$. To that purpose, it is
enough   to show    that $ \lim_{\epsilon\to  0}
|\lambda(\epsilon){\bfm u^\epsilon}|_\pi   =   |\overline{{\bfm
f}}|_\pi$. We multiply \eqref{multiplie-lambda} by
$\lambda(\epsilon)$ and  we deduce
$$
\overline{\lim}_{\epsilon\to  0} |\lambda(\epsilon){\bfm
u^\epsilon}|_\pi^2 \le \overline{\lim}_{\epsilon\to  0} ({\bfm
 f},  \lambda(\epsilon){\bfm u^\epsilon})_\pi    =  ({\bfm
 f},\overline{{\bfm f}})_\pi  =  |\overline{{\bfm f}}|_\pi^2
$$
which yields   the desired result   since  $|\overline{{\bfm
f}}|_\pi      \le   \underline{\lim}_{\epsilon\to  0}
|\lambda(\epsilon){\bfm u^\epsilon}|_\pi^2$ as  a consequence   of
the   weak convergence.\qed


\noindent {\bf Proof of Theorem \ref{ergth}.} It is enough to investigate the case of a function ${\bfm f}\in \mathcal{C}$. Indeed, the general case the results from the inequality
\begin{align*}
\M_\pi\E\Big[\sup_{0\leq t \leq T}\big|\int_0^t{\bfm f}(\tau_{\overline{X}^\epsilon_{r-}}\omega)\,dr-t\M_\pi[{\bfm f}]\big|\Big]\leq & \M_\pi\E\Big[\sup_{0\leq t \leq T}\big|\int_0^t{\bfm f}_n(\tau_{\overline{X}^\epsilon_{r-}}\omega)\,dr-t\M_\pi[{\bfm f}_n]\big|\Big]\\
&+2T\M_\pi[|{\bfm f}-{\bfm f}_n|].
\end{align*}
and the density of $\mathcal{C}$ in $L^1(\Omega)$.

So we consider a function ${\bfm f}\in \mathcal{C}$. Furthermore, even if it means replacing ${\bfm f}$ by ${\bfm f}-\M_\pi[{\bfm f}] $, we may (and will) assume that $\M_\pi[{\bfm f}]=0 $. Since $\mathcal{C}\subset H^\infty(\Omega)\cap L^\infty(\Omega)$, Proposition \ref{core} ensures that the solution ${\bfm u}^\epsilon$ of the resolvent equation $\lambda(\epsilon) {\bfm u}^\epsilon-{\bfm L}^\epsilon{\bfm u}^\epsilon = {\bfm f}$ belongs to $H^\infty(\Omega)\cap L^\infty(\Omega)$. So we can apply the It\^o
formula:
\begin{align*}
{\bfm  u}^\epsilon(\tau_{\overline{X}_{t}^\epsilon}\omega)=& {\bfm
u}^\epsilon(\tau_{x/\epsilon}\omega)  +   \int_0^t
\frac{1}{2\epsilon^2}{\bfm L^\epsilon}{\bfm
u}^\epsilon
(\tau_{\overline{X}_{r^-}^\epsilon}\omega)\,dr  + \int_0^t \frac{1}{\epsilon} {\bfm \sigma} D {\bfm  u}^\epsilon(\tau_{\overline{X}_{r^-}^\epsilon}\omega)\, d
 B_r \\& + \int_0^t \int_\R\Big[{\bfm  u}^\epsilon (\tau_{\overline{X}_{r-}^\epsilon +
 \gamma(\tau_{\overline{X}_{r-}^\epsilon}\omega, \frac{z}{\epsilon})}\omega)  -  {\bfm
 u}^\epsilon(\tau_{\overline{X}_{r-}^\epsilon}\omega)\Big] \widetilde{N}(d r ,
 dz)\nonumber.
\end{align*}
By using the relation $\lambda(\epsilon) {\bfm u}^\epsilon-{\bfm L}^\epsilon{\bfm u}^\epsilon = {\bfm f} $,  we  deduce    for  every  $\epsilon>0$
and $t\ge 0$,
\begin{align}
\int_0^t   {\bfm f}(\tau_{\overline{X}_{r^-}^\epsilon}\omega)  dr &=
-\epsilon^2 ({\bfm
u}^\epsilon(\tau_{\overline{X}_{t}^\epsilon}\omega) - {\bfm
u}^\epsilon(\tau_{x/\epsilon}\omega))  + \int_0^t
\lambda(\epsilon){\bfm
u}^\epsilon(\tau_{\overline{X}_{r-}^\epsilon}\omega) dr  +
\epsilon\int_0^t
{\bfm\sigma}  D {\bfm u}^\epsilon(\tau_{\overline{X}_{r-}^\epsilon}\omega)  dB_r\nonumber\\
&+\epsilon^2\int_0^t \int_\R\Big[{\bfm  u}^\epsilon
(\tau_{\overline{X}_{r^-}^\epsilon +
 \gamma(\tau_{\overline{X}_{r-}^\epsilon}\omega, \frac{z}{\epsilon})}\omega)  -  {\bfm
 u}^\epsilon(\tau_{\overline{X}_{r-}^\epsilon}\omega)\Big] \widetilde{N}(d r ,
 dz).\label{relation-f}
\end{align}
We now establish the convergence to $0$ of each term of the above right-hand side as $\epsilon\to 0$.
 From Proposition \ref{core}, we have
\begin{equation}\label{th:est1}
\M_\pi\big[\sup_{0\leq t \leq T}|\epsilon^2 {\bfm
u}^\epsilon(\tau_{\overline{X}_{r-}^\epsilon}\omega) |^2\big]\leq
\epsilon^4\lambda(\epsilon)^{-2}|{\bfm f}|_\infty^2.
\end{equation}
Furthermore, the Jensen inequality and Proposition \ref{prop:invar} yield
\begin{equation}\label{th:est2}
\M_\pi\big[\sup_{0\leq t \leq T}|\int_0^t
\lambda(\epsilon){\bfm
u}^\epsilon(\tau_{\overline{X}_{r-}^\epsilon}\omega) dr|^2\big]\leq T \lambda(\epsilon)^2|{\bfm u}^\epsilon|_\pi^2.
\end{equation}
Concerning the Brownian martingale, we  use in turn the Doob inequality, Proposition \ref{prop:invar}, Assumption \ref{ellipticity} and \eqref{est:erg} to obtain
\begin{align}
\M_\pi\E\Big[\sup_{0\leq t \leq T}\big|\epsilon\int_0^t {\bfm\sigma} D {\bfm
u}^\epsilon  (\tau_{\overline{X}_{r^-}^\epsilon}\omega) dB_r\big|^2\Big] \leq &4 T \epsilon^2({\bfm a}D{\bfm
u}^\epsilon,D{\bfm u}^\epsilon)_\pi\leq 4 TM_{\ref{ellipticity}}\epsilon^2|D{\bfm
u}^\epsilon|_\pi^2\nonumber\\
=&4 TM_{\ref{ellipticity}}(\epsilon^2/\lambda(\epsilon))\,\lambda(\epsilon)|D{\bfm
u}^\epsilon|_\pi^2\leq 4 TM_{\ref{ellipticity}}(\epsilon^2/\lambda(\epsilon))\,|{\bfm
f}|_\pi^2.\label{th:est3}
\end{align}
We treat the jump martingale with Lemma \ref{mart:poisson} and  \eqref{est:erg} 
\begin{align}
\M_\pi\E\Big[\sup_{0\leq t \leq T}\big| \epsilon^2\int_0^t \int_\R & \Big[{\bfm  u}^\epsilon
(\tau_{\overline{X}_{r^-}^\epsilon +
 \gamma(\tau_{\overline{X}_{r-}^\epsilon}\omega, \frac{z}{\epsilon})}\omega)  -  {\bfm
 u}^\epsilon(\tau_{\overline{X}_{r-}^\epsilon}\omega)\Big] \widetilde{N}(d r ,
 dz)\big|^2\Big] \nonumber\\&=4T\epsilon^{4-\alpha}\M\int_\R \big(T_z{\bfm
u}^\epsilon-{\bfm u}^\epsilon \big)^2{\bfm c}(\omega,z)\,\nu(dz)\nonumber\\
&\leq 4T\big(\epsilon^{2}/\lambda(\epsilon)\big)\,\,|{\bfm
f}|_\pi^2.\label{th:est4}
\end{align}
We choose now a family $(\lambda(\epsilon))_\epsilon$ such that
$$\lim_{\epsilon\to 0}\lambda(\epsilon)=0\quad \text{and}\quad \lim _{\epsilon\to 0}\epsilon^2/\lambda(\epsilon)=0.$$
Estimates \eqref{th:est1} \eqref{th:est2} \eqref{th:est3} \eqref{th:est4} and Proposition \ref{prop:erg} ensure that all the terms involved in the right-hand side of \eqref{relation-f} converge to $0$ in $L^1(\Omega;C([0,T]))$ as $\epsilon\to 0$. This completes the proof.\qed

\section{Construction of the correctors}\label{sec:correctors}

\begin{proposition}
For  any   $\epsilon>0$,  we  define  $\bfm  u_\epsilon$ as  the
solution of the  resolvent  equation $$ \epsilon^2 \bfm  u^\epsilon
- {\bfm L}^\epsilon{\bfm u}^\epsilon   = {\bfm  b},
$$
that is $ \bfm  u^\epsilon=G^\epsilon_{\epsilon^2}\bfm b$. Then we
can find $\bfm  \xi   \in L^2(\Omega)$ such that
$$ \lim_{\epsilon\to 0}\Big[ \epsilon^2 |{\bfm u}^\epsilon|_\pi^2 + |D
{\bfm u}^\epsilon  -   {\bfm \xi}|_\pi^2  + \epsilon^{2-\alpha}\M
\int_\R (T_z {\bfm  u}^\epsilon - {\bfm u}^\epsilon)^2 {\bfm
c}(\omega, z)\nu(dz)\Big] =   0
$$
\end{proposition}
\noindent {\it Proof.} By plugging the function $\bfm \psi=\bfm  u^\epsilon
$ into  the weak form of the resolvent equation \eqref{eqbase}, we
get
\begin{align}
\epsilon^2 |\bfm  u^\epsilon|_\pi^2  +  \frac{1}{2}({\bfm  a} D
{\bfm u^\epsilon}, D {\bfm  u^\epsilon})_\pi  &+
\frac{\epsilon^{2-\alpha}}{2}\M\int_\R (T_z {\bfm u}^\epsilon -
{\bfm u}^\epsilon)^2 {\bfm c}(\omega, z)\nu(dz) =
({\bfm  b},{\bfm u}^\epsilon)_\pi \label{multiplie-u}
%
\end{align}
We estimate the right-hand side by
\begin{align*}
({\bfm  b},{\bfm u}^\epsilon)_\pi =\frac{1}{2}(D({\bfm  a}e^{-2{\bfm
V}}),{\bfm u}^\epsilon)_\pi=-\frac{1}{2}({\bfm  a}e^{-2{\bfm
V}},D{\bfm u}^\epsilon)_\pi=-\frac{1}{2}({\bfm  a},D{\bfm
u}^\epsilon)_\pi.
\end{align*}
By using the Cauchy Schwarz inequality, we obtain for some constant
$C>0$ (independent of $\epsilon $)
$$ ({\bfm  b},{\bfm u}^\epsilon)_\pi\leq C+\frac{1}{4}({\bfm  a} D {\bfm
u^\epsilon}, D {\bfm  u^\epsilon})_\pi.$$ Plugging this inequality
into \eqref{multiplie-u} yields
$$
\epsilon^2 |\bfm  u^\epsilon|_\pi^2  +  \frac{1}{4}({\bfm  a} D
{\bfm u^\epsilon}, D {\bfm  u^\epsilon})_\pi  +
\frac{\epsilon^{2-\alpha}}{2}\M\int_\R (T_z {\bfm u}^\epsilon -
{\bfm u}^\epsilon)^2 {\bfm c}(\omega, z)\nu(dz)  \le C
$$
Assumption \ref{ellipticity} then ensures that the family $(D
{\bfm u^\epsilon})_\epsilon$ is bounded in $L^2(\Omega) $. Therefore   we can find $ {\bfm \xi}\in L^2(\Omega)$ such that,
along a  subsequence, the family $(D {\bfm
u^\epsilon})_{\epsilon>0}$ converges weakly in $L^2(\Omega)$ towards
$ {\bfm \xi}$ as $\epsilon\to 0$. We further stress that the
previous bound implies that
$$
\epsilon^2 \bfm  u^\epsilon \xrightarrow[L^2(\Omega)]{strongly}  0;
\quad  \quad    \epsilon^{2-\alpha} (T_z {\bfm u}^\epsilon - {\bfm
u}^\epsilon)\xrightarrow[L^2(\Omega\times\R;{\bfm
c}(\omega,z)\nu(dz)\mu(d\omega))]{strongly} 0.
$$

Now we establish that the whole family $(D {\bfm
u^\epsilon})_{\epsilon>0}$ is weakly converging. From the   resolvent equation \eqref{eqbase}, we  have for any
$\varphi \in \H$,
\begin{equation}\label{eq:rescor}
(\epsilon^2 \bfm  u^\epsilon, \bfm \varphi)_\pi   +
\frac{1}{2}({\bfm a} D {\bfm u^\epsilon}, D \bfm\varphi)_\pi  +
\frac{\epsilon^{2-\alpha}}{2}\M\int_\R (T_z {\bfm u}^\epsilon -
{\bfm u}^\epsilon) (T_z \bfm\varphi - \bfm\varphi) {\bfm c}(\omega,
z)\nu(dz) = ({\bfm  b},\bfm\varphi)_\pi
\end{equation}
Letting  $\epsilon\to  0$ yields for any $\varphi \in  \H, \quad
\frac{1}{2}({\bfm a}  {\bfm \xi}, D \bfm \varphi)_\pi=({\bfm b},\bfm
\varphi)_\pi$.  This equation characterizes the function ${\bfm \xi}$ and therefore establishes the uniqueness of the weak limit.
So the whole family $(D {\bfm u^\epsilon})_{\epsilon>0}$ is weakly converging.

We now establish the strong convergence. As a consequence of the
previous equality, we have $$\lim_{\epsilon\to 0}({\bfm b},{\bfm
u}^\epsilon)_\pi=\lim_{\epsilon\to 0}\frac{1}{2}({\bfm a} {\bfm
\xi}, D {\bfm u}^\epsilon)_\pi=\frac{1}{2}({\bfm a}  {\bfm \xi},
{\bfm \xi})_\pi.$$ We take the $\overline{\lim}_{\epsilon\to 0}$ in
both sides of \eqref{multiplie-u}, we obtain
\begin{equation}\label{eq:limsup}
\overline{\lim}_{\epsilon\to  0} \Big[\epsilon^2 |\bfm
u^\epsilon|_\pi^2  +  \frac{1}{2}({\bfm  a} D {\bfm u^\epsilon}, D
{\bfm  u^\epsilon})_\pi  + \frac{\epsilon^{2-\alpha}}{2}\M\int_\R
(T_z {\bfm u}^\epsilon - {\bfm u}^\epsilon)^2 {\bfm c}(\omega,
z)\nu(dz)\Big] \le \frac{1}{2}({\bfm a}  {\bfm \xi}, {\bfm \xi})_\pi.
\end{equation}
We deduce
\begin{equation}\label{eq:limsup1}
\overline{\lim}_{\epsilon\to 0}({\bfm  a} D {\bfm u^\epsilon}, D
{\bfm  u^\epsilon})_\pi \leq ({\bfm a}  {\bfm \xi}, {\bfm \xi})_\pi.
\end{equation}
Note that both inner products $(\cdot,\cdot)_\pi $ and $({\bfm
a}\cdot,\cdot)_\pi $ define equivalent norms on $L^2(\Omega) $ since
$\bfm a$ is uniformly elliptic (Assumption \ref{ellipticity}). The
family $(D{\bfm  u^\epsilon})_{\epsilon>0} $ therefore weakly
converges towards $ \bfm \xi$ in $L^2(\Omega) $ with respect to the
norm associated to $({\bfm a}\cdot,\cdot)_\pi $. The weak
convergence implies \begin{equation}\label{eq:limsup2}
\underline{\lim}_{\epsilon\to 0}({\bfm a} D {\bfm u^\epsilon}, D
{\bfm  u^\epsilon})_\pi \geq ({\bfm a} {\bfm \xi}, {\bfm \xi})_\pi.
\end{equation}
Gathering \eqref{eq:limsup1} and \eqref{eq:limsup2} yields
$\lim_{\epsilon\to 0}({\bfm a} D {\bfm u^\epsilon}, D {\bfm
u^\epsilon})_\pi = ({\bfm a} {\bfm \xi}, {\bfm \xi})_\pi $. The
convergence of the norms together with the weak convergence implies
the strong convergence of the family $(D{\bfm
u^\epsilon})_{\epsilon>0} $  towards $ \bfm \xi$ in $L^2(\Omega) $
with respect to the norm associated to $({\bfm a}\cdot,\cdot)_\pi $.
Since the norm associated to $({\bfm a}\cdot,\cdot)_\pi $ is
equivalent to that associated to $(\cdot,\cdot)_\pi $, we deduce the
strong convergence of the family $(D{\bfm  u^\epsilon})_{\epsilon>0}
$ towards $ \bfm \xi$ in $L^2(\Omega) $. Once that convergence
established, \eqref{eq:limsup} also implies
$$\overline{\lim}_{\epsilon\to  0} \Big[\epsilon^2 |\bfm
u^\epsilon|_\pi^2  + \frac{\epsilon^{2-\alpha}}{2}\M\int_\R (T_z
{\bfm u}^\epsilon - {\bfm u}^\epsilon)^2 {\bfm c}(\omega,
z)\nu(dz)\Big] =0,$$ and this completes the proof. \qed

\section{Homogenization}\label{sec:homo}


We apply the It\^o formula to the function $x\mapsto \epsilon{\bfm u}_\epsilon(\tau_{x/\epsilon}\omega)$ where ${\bfm u}_\epsilon=G^\epsilon_{\epsilon^2}({\bfm b})$ (see Section \ref{sec:correctors}) and we get:
\begin{align}
\epsilon {\bfm u}^\epsilon(\tau_{\overline{X}^\epsilon_t}\omega)-\epsilon {\bfm u}^\epsilon(\tau_{x/\epsilon}\omega)=&\int_0^t \epsilon {\bfm u}^\epsilon(\tau_{\overline{X}^\epsilon_{r-}}\omega) dr
-\int_0^t \frac{1}{\epsilon}{\bfm b}(\tau_{\overline{X}^\epsilon_{r-}}\omega) dr+\int_0^tD{\bfm u}^\epsilon{\bfm \sigma}(\tau_{\overline{X}^\epsilon_r}\omega)\,dB_r \nonumber\\
\label{itou}&+\epsilon\int_0^t\int_\R\big({\bfm
u}^\epsilon(\tau_{\overline{X}^\epsilon_{r-}+{\bfm
\gamma}(\tau_{\overline{X}^\epsilon_{r-}}\omega,z)}\omega)-{\bfm
u}^\epsilon
(\tau_{\overline{X}^\epsilon_{r-}}\omega)\big)\,\widetilde{N}(dr,dz).
\end{align}
Therefore, by summing with \eqref{SDE}, we get:
\begin{align*}
\epsilon{\bfm u}^\epsilon&(\tau_{\overline{X}^\epsilon_t}\omega)+X^\epsilon_t=x+\epsilon{\bfm u}^\epsilon(\tau_{x/\epsilon}\omega)
+\int_0^t \epsilon {\bfm u}^\epsilon(\tau_{\overline{X}^\epsilon_{r-}}\omega) dr +\int_0^t(1+D{\bfm u}^\epsilon){\bfm \sigma}(\tau_{\overline{X}^\epsilon_{r-}}\omega)\,dB_r\\
&+\epsilon\int_0^t\int_\R\big({\bfm u}^\epsilon(\tau_{\overline{X}^\epsilon_{r-}+{\bfm \gamma}(\tau_{\overline{X}^\epsilon_{r-}}\omega,\frac{z}{\epsilon})}\omega)-
{\bfm u}^\epsilon (\tau_{\overline{X}^\epsilon_{r-}}\omega)\big)\,\widetilde{N}(dr,dz)\\
&+\epsilon\int_0^t\int_\R{\bfm
\gamma}(\tau_{\overline{X}^\epsilon_{r-}}\omega,\frac{z}{\epsilon})\,\hat{N}^\epsilon(dr,dz)+\int_0^t
\epsilon^{1-\alpha} {\bfm
e}(\tau_{\overline{X}^\epsilon_{r-}}\omega) dr
\end{align*}
In order to prove the result, we consider each term in the above sum separately.

\begin{lemma}\label{cvzero}
We have the following convergence:
$$\M_\pi\big[|\epsilon{\bfm u}^\epsilon(\tau_{\overline{X}^\epsilon_t}\omega)|^2+|\epsilon{\bfm u}^\epsilon(\tau_{x/\epsilon}\omega)|^2
+\int_0^t \epsilon {\bfm
u}^\epsilon(\tau_{\overline{X}^\epsilon_{r-}}\omega) dr|^2\big]\to
0, \quad \text{as }\epsilon\to 0$$ and
$$\M_\pi\E\Big[|\epsilon\int_0^t\int_\R \big({\bfm u}^\epsilon(\tau_{\overline{X}^\epsilon_{r-}+{\bfm \gamma}(\tau_{\overline{X}^\epsilon_{r-}}\omega,
\frac{z}{\epsilon})}\omega)-{\bfm u}^\epsilon
(\tau_{\overline{X}^\epsilon_{r-}}\omega)\big)\,\widetilde{N}(dr,dz)|^2\Big]
\to 0, \quad \text{as }\epsilon\to 0. $$
\end{lemma}

Thus, we just have to investigate the convergence of the following semimartingale:
\begin{equation*}
Y^{\epsilon}_t=\int_0^t(1+D{\bfm u}^\epsilon){\bfm
\sigma}(\tau_{\overline{X}^\epsilon_{r-}}\omega)\,dB_r+\epsilon\int_0^t\int_\R{\bfm
\gamma}(\tau_{\overline{X}^\epsilon_{r-}}\omega,\frac{z}{\epsilon})\,\hat{N}^\epsilon(dr,dz)+\int_0^t
\epsilon^{1-\alpha} {\bfm
e}(\tau_{\overline{X}^\epsilon_{r-}}\omega) dr
\end{equation*}
In order to obtain the desired result, we introduce the truncation function $h$
defined by $h(z)=z$ if $|z|\leq 1$ and $h(z)={\rm sign}(z)$ if $|z|>1$ and we use theorem VIII.4.1 in \cite{jacod}.
Following the notations of \cite{jacod}, we introduce the following processes:
\begin{equation*}
\check{Y}^{\epsilon,(h)}_t=\sum_{0<s \leq t}\Delta Y^{\epsilon}_s-h(\Delta Y^{\epsilon}_s)
\end{equation*}
and
\begin{equation*}
Y^{\epsilon,(h)}_t=Y^{\epsilon}_t-\check{Y}^{\epsilon,(h)}_t.
\end{equation*}
Note that we can decompose the semimartingale $Y^{\epsilon,(h)}$ into its martingale part and its predictable bounded variation part as:
\begin{equation*}
Y^{\epsilon}_t=M^{\epsilon,(h)}_t+B^{\epsilon,(h)}_t,
\end{equation*}
where $M^{\epsilon,(h)},B^{\epsilon,(h)}$ are given by:
\begin{equation*}
M^{\epsilon,(h)}_t= \int_0^t(1+D{\bfm u}^\epsilon){\bfm
\sigma}(\tau_{\overline{X}^\epsilon_{r-}}\omega)\,dB_r+\int_0^{ t }
\int_{\R}
h(\epsilon\gamma(\tau_{\overline{X}_{r-}}\omega,\frac{z}{\epsilon}))\,\widetilde{N}(dr,dz)
\end{equation*}
and
\begin{align*}
B^{\epsilon,(h)}_t= &\int_0^t \epsilon^{1-\alpha} {\bfm
e}(\tau_{\overline{X}^\epsilon_{r-}}\omega) dr +\int_0^{ t }
\int_{\R}
h(\epsilon\gamma(\tau_{\overline{X}_{r-}}\omega,\frac{z}{\epsilon}))\one_{\{|z|>\epsilon\}}\,\nu(dz)dr
\end{align*}
According to  theorem VIII.4.1 in \cite{jacod}, to prove the convergence of the semimartingale $Y^\epsilon$ in the Skorohod topology towards a L\'evy process with characteristic function given by Theorem \ref{mainresult}, we have to establish:

1) $ B^{\epsilon,(h)}$  converges towards $0$ in $C([0,T];\R)$ with respect to the sup-norm in probability,

2)  we denote by $<M^{\epsilon,(h)}>$ the compensator of the martingale $M^{\epsilon,(h)}$, that is the unique $\mathcal{F}_t$-predictable process such that $M^{\epsilon,(h)}-<M^{\epsilon,(h)}>$ becomes a $\mathcal{F}_t$-martingale. Then $<M^{\epsilon,(h)}>$ converges towards $At+t\M[{\bfm \theta}]\int_\R h^2(z)\nu(dz) $ in probability for all $t\geq 0$.

3) for every bounded continuous function $g$ vanishing in a neighborhood of $0$, the following convergence holds:
$$\int_0^t\int_\R g\big(\epsilon{\bfm \gamma}(\tau_{\overline{X}^\epsilon_{r-}}\omega,\frac{z}{\epsilon})\big)\nu(dz)dr\to t\M[{\bfm \theta}]\int_\R g(z)\nu(dz)$$
in probability for all $t\geq 0$.

So we check now that all the above points are satisfied. Point 3) clearly results from corollary \ref{corerg}. The compensator of the martingale $M^{\epsilon,(h)}$ is given by
\begin{align*}
<M^{\epsilon,(h)}>_t= &\int_0^t(1+D{\bfm u}^\epsilon)^2{\bfm
a}(\tau_{\overline{X}^\epsilon_{r-}}\omega)\,dr+\int_0^{ t }
\int_{\R}
h^2(\epsilon\gamma(\tau_{\overline{X}_{r-}}\omega,\frac{z}{\epsilon}))\,\nu(dz)dr.
\end{align*}
Point 2) then results from the combination of Section
\ref{sec:correctors} (in particular $D{\bfm u}^\epsilon\to {\bfm
\xi}$ in $L^2(\Omega)$, Corollary \ref{corerg2} and Corollary
\ref{corerg}.

We admit for a while the following result, which is proved in the appendix.
\begin{lemma}\label{lem:b}
The process $B^{\epsilon,(h)} $ converges to $0$ in $D([0,T];\R)$ in probability for the Skorohod topology.
\end{lemma}

To sum up, the three characteristics of the semimartingale $Y^\epsilon$ converge as $\epsilon\to 0$ to those of a L\'evy process $L$ with L\'evy exponent:
\begin{equation*}
\varphi(u)=-\frac{1}{2}Au^2+\int_{\R}(e^{iuz}-1-iuz\one_{\{|z|\leq 1\}})\M[{\bfm \theta}]\nu(dz).
\end{equation*}
Using theorem VIII.4.1 in \cite{jacod}, we conclude that the following convergence holds for the Skorohod topology:
\begin{equation*}
Y^{\epsilon} \overset{\epsilon \to 0}{\longrightarrow} L.
\end{equation*}
We deduce that the finite-dimensional distributions of the process $X^\epsilon$ converge to those of the process $L$. It remains to prove that the process $X^\epsilon$ is tight for the Skorohod topology. This is the purpose of the next section.\qed


\section{Tightness}\label{tightness}


Arguing exactly as in Section \ref{sec:homo}, we can prove that the semimartingale
\begin{equation*}
Y^{\epsilon}_t=\int_0^t{\bfm
\sigma}(\tau_{\overline{X}^\epsilon_{r-}}\omega)\,dB_r+\epsilon\int_0^t\int_\R{\bfm
\gamma}(\tau_{\overline{X}^\epsilon_{r-}}\omega,\frac{z}{\epsilon})\,\hat{N}^\epsilon(dr,dz)+\int_0^t
\epsilon^{1-\alpha} {\bfm
e}(\tau_{\overline{X}^\epsilon_{r-}}\omega) dr
\end{equation*}
converges for the Skorohod topology as $\epsilon\to 0$. So it is tight. The difficult term actually is
$$\frac{1}{\epsilon}\int_0^t{\bfm b}(\tau_{\overline{X}^\epsilon_{r-}}\omega)\,dr. $$
The strategy to establish its tightness is inspired from \cite[Section 3.3]{olla} (idea originally adapted from \cite{varadhan}).  The adaptation to the setup of jump-diffusion processes is given in \cite{vargas}. The present setup does not give rise to additional difficulties. So we just outline the proof and write properly the intermediate steps to stick with the notations of the present paper. The reader is referred to \cite{vargas} or \cite{olla} for further details.

Our purpose is to establish the following result:
\begin{theorem}
We have  the following estimation of the continuity modulus:
\begin{equation}
\M_\pi\E\Big[\sup_{\substack{|t-s|\leq \delta\\0\leq s,t\leq
T}}\big|\frac{1}{\epsilon}\int_0^t{\bfm
b}(\tau_{\overline{X}^\epsilon_{r-}}\omega)\,dr\big|\Big]\leq
C(T)\delta^{1/2}\ln \delta^{-1}
\end{equation}
for some positive constant $C(T)$ only depending on $T$.
\end{theorem}
The tightness of $\frac{1}{\epsilon}\int_0^t{\bfm
b}(\tau_{\overline{X}^\epsilon_{r-}}\omega)\,dr$ in the Skorohod
topology is a direct consequence.

\vspace{2mm}
\noindent {\bf Guideline of the proof.} We assume that the starting point $x$ is $0$ in \eqref{SDE}. We will explain thereafter how to deduce the general case.

1) Since ${\bfm b}=\frac{e^{2\bfm V}}{2}D({\bfm a}e^{-2\bfm V})$, we can make an integration by parts to check that
\begin{equation}\label{poincaretight}
\forall {\bfm \varphi}\in \mathcal{C},\quad ({\bfm b},{\bfm \varphi}^2)_\pi\leq P B^d( {\bfm \varphi}, {\bfm \varphi})^{1/2}| {\bfm \varphi}|_\pi \leq P \big(B^d( {\bfm \varphi}, {\bfm \varphi})+\epsilon^{2-\alpha}B^j( {\bfm \varphi}, {\bfm \varphi}) \big)^{1/2}| {\bfm \varphi}|_\pi
\end{equation}
for some positive constant $P$.

2) Then we estimate the exponential moments of the random variable
$\frac{1}{\epsilon}\int_0^t{\bfm
b}(\tau_{\overline{X}^\epsilon_{r-}}\omega)\,dr$. To that purpose,
the Feynmann-Kac formula provides a connection between the
exponential moments and the solution of a certain evolution
equation:
\begin{theorem}{\bf Feynmann-Kac formula.}\label{fk}
Fix $ \beta\in\R$. The function
$${\bfm u}^\epsilon(t,\omega)=\E\Big[\exp\big(\frac{\beta}{\epsilon}\int_0^t{\bfm b}(\tau_{\overline{X}^\epsilon_{r-}}\omega)\,dr\big)\Big]$$ is a solution of the equation
$$ \partial_t{\bfm u}^\epsilon=\epsilon^{-2}{\bfm L}^\epsilon{\bfm u}^\epsilon+\beta\epsilon^{-1}{\bfm b}{\bfm u}^\epsilon$$
with initial condition ${\bfm u}^\epsilon(0,\omega)={\bfm 1}$.
\end{theorem}
\begin{remark}
By solution, we mean a function ${\bfm u}^\epsilon$ such that
$\forall t \geq 0$, ${\bfm u}^\epsilon(t,\cdot)\in{\rm Dom}({\bfm
L}^\epsilon)$ and $$\lim_{s\to 0}\frac{{\bfm
u}^\epsilon(t+s,\cdot)-{\bfm
u}^\epsilon(t,\cdot)}{s}=\epsilon^{-2}{\bfm L}^\epsilon{\bfm
u}^\epsilon(t,\cdot)+\beta\epsilon^{-1}{\bfm b}(\cdot){\bfm
u}^\epsilon(t,\cdot)\quad \text{in }L^2(\Omega).$$
\end{remark}
By using the Dirichlet form associated to the operator ${\bfm L}^\epsilon$, we can prove
\begin{proposition}\label{carreu}
Let ${\bfm u}^\epsilon(t,\cdot)$ be the function of Theorem
\ref{fk}. Then
$$ \M_\pi[{\bfm u}^\epsilon(t,\cdot)^2]\leq e^{2\lambda_\epsilon t}$$ where $\lambda_\epsilon$ is defined as  $\lambda_\epsilon=\sup_{\substack{|{\bfm \varphi}|_\pi=1,\\{\bfm \varphi}\in{\rm Dom}{\bfm L}^\epsilon}}({\bfm \varphi},(\epsilon^{-2}{\bfm L}^\epsilon+\beta\epsilon^{-1}{\bfm b}){\bfm \varphi})_\pi$.
\end{proposition}

By using the stationnarity of the process
$\tau_{\overline{X}^\epsilon_r} \omega$ under the measure $\pi$ and
Proposition \ref{carreu}, we then establish
 \begin{equation}\label{estbeta}
 \M_\pi\E\Big[\exp\Big|\beta\frac{1}{\epsilon}\int_{s}^{t}{\bfm b}(\tau_{\overline{X}^\epsilon_{r-}}\omega)\,dr\Big|\Big] \leq 2\exp\Big(\lambda_\epsilon(t-s)\Big).
 \end{equation}
By using \eqref{poincaretight}, we get for each function ${\bfm \varphi}\in {\rm Dom}({\bfm L}^\epsilon) $ such that $|{\bfm \varphi}|_\pi=1$:
\begin{align*}
({\bfm \varphi},(\epsilon^{-2}{\bfm L}^\epsilon+\beta\epsilon^{-1}{\bfm b}){\bfm \varphi})_\pi \leq &-\epsilon^{-2}\big(B^d({\bfm \varphi},{\bfm \varphi})+\epsilon^{2-\alpha}B^j({\bfm \varphi},{\bfm \varphi})\big)\\&+P\beta\epsilon^{-1} \big(B^d({\bfm \varphi},{\bfm \varphi})+\epsilon^{2-\alpha}B^j({\bfm \varphi},{\bfm \varphi})\big)^{1/2}.
\end{align*}
Optimizing the polynom $P(x)=-\epsilon^{-2}x^2+P\beta \epsilon^{-1} x $ with respect to the variable $x\in\R$ yields
\begin{equation}\label{estl}
\lambda_\epsilon \leq \beta^2P^2/4.
\end{equation}
We gather \eqref{estbeta} and \eqref{estl} to obtain
$$\M_\pi\E\Big[\exp\Big|\beta\frac{1}{\epsilon}\int_{s}^{t}{\bfm b}(\tau_{\overline{X}^\epsilon_{r-}}\omega)\,dr\Big|\Big] \leq 2 \exp\Big(\beta^2P^2(t-s)/4\Big).$$

3) The last step is to use the GRR inequality to exploit the exponential bounds we have just established
\begin{proposition}\label{rumsey}{\bf (Garsia-Rodemich-Rumsey's
inequality).} Let $p$ and $\Psi$ be strictly increasing continuous
functions on $[0,+\infty[$ satisfying $p(0)=\Psi(0)=0$ and
$\lim_{t\to \infty}\Psi(t)=+\infty$. For given $T>0$ and $g\in
C([0,T];\R^d)$, suppose that there exists a finite $B$ such that;
\begin{equation}\label{GRRcond}
  \int_0^T\int_0^T\Psi\Big(\frac{|g(t)-g(s)|}{p(|t-s|)}\Big)\,ds\,dt\leq B<\infty.
\end{equation}
Then, for all $0\leq s \leq t \leq T$,
\begin{equation}\label{GRR}
  |g(t)-g(s)|\leq 8\int_0^{t-s}\Psi^{-1}(4B/u^2)\,dp(u).
\end{equation}
\end{proposition}

We conclude by using the GRR inequality (with
$g(t)=\frac{1}{\epsilon}\int_{s}^{t}{\bfm
b}(\tau_{\overline{X}^\epsilon_{r-}}\omega)\,dr$, $p(t)=\sqrt{t}$,
$\Psi(t)=e^t-1$), by taking the expectation and by using the above
estimate.

4) In the case the starting point $x$ is not necessary equal to $0$. We denote by $X^{\epsilon,\omega,x}$ the solution of \eqref{SDE} starting from $x\in\R$ in the environment $\omega$. It is plain to see that the processes $x+X^{\epsilon,\tau_{\frac{x}{\epsilon}}\omega,0}$ and $X^{\epsilon,\omega,x }$ are both solution of the same SDE. So they have the same law. In what follows, $C$ denotes a positive constant such that $C^{-1}\M[{\bfm f}]\leq \M_\pi[{\bfm f}] \leq C\M[{\bfm f}] $ for all measurable nonnegative functions ${\bfm f}:\Omega\to \R $. We deduce
\begin{align*}
\M_\pi\E\Big[\sup_{\substack{|t-s|\leq \delta\\0\leq s,t\leq
T}}\big|\frac{1}{\epsilon}\int_0^t{\bfm
b}(\tau_{\overline{X}^{\epsilon,\omega,x}_{r-}}\omega)\,dr\big|\Big]&=\M_\pi\E\Big[\sup_{\substack{|t-s|\leq
\delta\\0\leq s,t\leq T}}\big|\frac{1}{\epsilon}\int_0^t{\bfm
b}(\tau_{
\frac{x}{\epsilon}+\overline{X}^{\epsilon,\tau_{x/\epsilon}\omega,0}_{r-}
}\,\omega)\,dr\big|\Big] \\
&=\M_\pi\E\Big[\sup_{\substack{|t-s|\leq \delta\\0\leq s,t\leq
T}}\big|\frac{1}{\epsilon}\int_0^t{\bfm b}\big(\tau_{
\overline{X}^{\epsilon,\tau_{x/\epsilon}\omega,0}_{r-}
}\,(\tau_{\frac{x}{\epsilon}}\omega)\big)\,dr\big|\Big] \\
&\leq C\M\E\Big[\sup_{\substack{|t-s|\leq \delta\\0\leq s,t\leq
T}}\big|\frac{1}{\epsilon}\int_0^t{\bfm b}\big(\tau_{
\overline{X}^{\epsilon,\tau_{x/\epsilon}\omega,0}_{r-}
}\,(\tau_{\frac{x}{\epsilon}}\omega)\big)\,dr\big|\Big]
\end{align*}
By invariance of the measure $\mu$ under translations, the last quantity matches
$$ C\M\E\Big[\sup_{\substack{|t-s|\leq \delta\\0\leq s,t\leq T}}\big|\frac{1}{\epsilon}\int_0^t{\bfm b}\big(\tau_{
\overline{X}^{\epsilon,\omega,0}_{r-} }\,\omega\big)\,dr\big|\Big]
.$$ Finally, we deduce
\begin{align*}
\M_\pi\E\Big[\sup_{\substack{|t-s|\leq \delta\\0\leq s,t\leq
T}}\big|\frac{1}{\epsilon}\int_0^t{\bfm
b}(\tau_{\overline{X}^{\epsilon,\omega,x}_{r-}}\omega)\,dr\big|\Big]\leq
&   C^2\M_\pi\E\Big[\sup_{\substack{|t-s|\leq \delta\\0\leq s,t\leq
T}}\big|\frac{1}{\epsilon}\int_0^t{\bfm b}\big(\tau_{
\overline{X}^{\epsilon,\omega,0}_{r-} }\,\omega\big)\,dr\big|\Big]\\
\leq & C(T)\delta^{1/2}\ln \delta^{-1}.
\end{align*}
So we have established the tightness estimates for any starting point $x\in\R$.
\qed


\appendix
\section*{Appendix}
\section{Auxiliary lemmas}


\begin{lemma}\label{limitz}
If the function $h:\R\to \R$ is given by $h(z)=z$ if $|z|\leq 1$ and $h(z)={\rm sign}(z)$ otherwise, then the following limit holds in the  $L^\infty(\Omega)$ sense
$$\lim_{\beta\downarrow 0}\int_{\beta<|z|}h(z){\bfm c}(\omega,\frac{z}{\epsilon})e^{2{\bfm V}(\omega)}\nu(dz)$$ uniformly with respect to $ \epsilon>0$.
\end{lemma}

\noindent {\it Proof.} It it plain to adapt the proof of Lemma A.2 in \cite{vargas} with ${\bfm g}(\omega,z)=h(z)e^{2{\bfm V}(\omega)}$ and $k=0$. The function  ${\bfm C}$ is here given by ${\bfm C}(\omega)=2{\bfm V}(\omega) $ since $h(z)+h(-z)=0$. Following the proof in  \cite{vargas}, it is not difficult to see that the convergence is uniform with respect to $ \epsilon$. We let the reader check the details.\qed

\begin{lemma}\label{mart:poisson}
If ${\bfm \varphi}\in\H $ then
\begin{align*}
\M_\pi\Big[\sup_{0\leq t \leq T}\big|\int_0^t\int_\R  \big({\bfm \varphi}(\tau_{\overline{X}^\epsilon_{r-}+{\bfm \gamma}(\tau_{\overline{X}^\epsilon_{r-}}\omega,\frac{z}{\epsilon})}\omega) & -{\bfm \varphi} (\tau_{\overline{X}^\epsilon_{r-}}\omega)\big)\widetilde{N}(dr,dz)\big|^2\Big] \\
&\leq 4T \epsilon^{-\alpha}\int_0^t\M \big[\int_\R \big(T_y{\bfm
\varphi}-{\bfm \varphi} \big)^2{\bfm c}(\omega,y)\,\nu(dy)\big].
\end{align*}
\end{lemma}

\noindent {\it Proof.} By using in turn the Doob inequality, Proposition \ref{prop:invar},
the change of variables $ z/\epsilon=y$ and the relation $\nu\circ {\bfm \gamma}_\omega^{-1}={\bfm c}(\omega,z)e^{2{\bfm V}(\omega)}\nu(dz)$, we obtain
\begin{align*}
\M_\pi\Big[\sup_{0\leq t \leq T}\big|\int_0^t\int_\R  \big({\bfm \varphi}(\tau_{\overline{X}^\epsilon_{r-}+{\bfm \gamma}(\tau_{\overline{X}^\epsilon_{r-}}\omega,\frac{z}{\epsilon})}\omega)&-{\bfm \varphi} (\tau_{\overline{X}^\epsilon_{r-}}\omega)\big)\widetilde{N}(dr,dz)\big|^2\Big]\\
& \leq 4\int_0^T\int_\R  \M_\pi\Big[\big|{\bfm \varphi}(\tau_{{\bfm \gamma}(\omega,\frac{z}{\epsilon})}\omega)-{\bfm \varphi} (\omega)\big|^2\Big]\,\nu(dz)\\
&=4T \epsilon^{-\alpha}\int_0^t\M \big[\int_\R \big(T_y{\bfm
\varphi}-{\bfm \varphi} \big)^2{\bfm c}(\omega,y)\,\nu(dy)\big]\qed
\end{align*}




\section{Proofs}

\noindent {\it Proof of Lemma \ref{cvzero}.} By using Section
\ref{measure} and the properties of the correctors (see Section
\ref{sec:correctors}), we have
$$\M_\pi\E\big[|\epsilon{\bfm u}^\epsilon(\tau_{\overline{X}^\epsilon_t}\omega)|^2+|\epsilon{\bfm u}^\epsilon(\tau_{x/\epsilon}\omega)|^2\big]\leq
2\epsilon^2|{\bfm u}^\epsilon|^2_\pi\to 0, \quad \text{as
}\epsilon\to 0.$$ Similarly,
$$\M_\pi\E\big[|\int_0^t \epsilon {\bfm u}^\epsilon(\tau_{\overline{X}^\epsilon_{r-}}\omega) dr|^2\big]\leq
 t\int_0^t \M_\pi\E\big[|\epsilon {\bfm u}^\epsilon(\tau_{\overline{X}^\epsilon_{r-}}\omega)|^2\big]  dr \leq
 t\epsilon^2|{\bfm u}^\epsilon|^2_\pi\to 0, \quad \text{as }\epsilon\to 0.$$
Concerning the jump martingale, we use Lemma \ref{mart:poisson}
\begin{align*}
\M_\pi\E\Big[|\epsilon\int_0^t\int_\R & \big({\bfm u}^\epsilon(\tau_{\overline{X}^\epsilon_{r-}+{\bfm \gamma}(\tau_{\overline{X}^\epsilon_{r-}}\omega,
\frac{z}{\epsilon})}\omega)-{\bfm u}^\epsilon (\tau_{\overline{X}^\epsilon_{r-}}\omega)\big)\,\widetilde{N}(dr,dz)|^2\Big] \\
& \leq \epsilon^{2-\alpha}t\M\int_\R  \big(T_z{\bfm
u}^\epsilon-{\bfm u}^\epsilon \big)^2{\bfm c}(\omega,z)\,\nu(dz)\to
0, \quad \text{as }\epsilon\to 0
\end{align*}\qed

\noindent {\it Proof of Lemma \ref{lem:b}.} The expression of
$B^{\epsilon,(h)}$ can be rewritten as $\int_0^t  {\bfm
g}_\epsilon(\tau_{\overline{X}^\epsilon_{r-}}\omega) dr$ where
\begin{align*}
{\bfm g}_\epsilon(\omega)= & \epsilon^{1-\alpha} {\bfm e}(\omega)  +\int_{\R}  h(\epsilon\gamma(\omega,\frac{z}{\epsilon}))\one_{\{|z|>\epsilon\}}\,\nu(dz)\\
= & \lim_{\beta\to 0}\int_{\{\beta\leq |{\bfm \gamma}(\omega,\frac{z}{\epsilon})|\}}\epsilon{\bfm \gamma}(\omega,\frac{z}{\epsilon})\one_{\{|z|\leq \epsilon\}}\nu(dz)  +\int_{\R}  h(\epsilon\gamma(\omega,\frac{z}{\epsilon}))\one_{\{|z|>\epsilon\}}\,\nu(dz)\\
= & \lim_{\beta\to 0}\int_{\{\beta\leq |{\bfm \gamma}(\omega,\frac{z}{\epsilon})|\}}  h(\epsilon\gamma(\omega,\frac{z}{\epsilon}))\,\nu(dz),
\end{align*}
the last equality resulting from $\sup_{|z|\leq \epsilon}|\epsilon{\bfm \gamma}(\omega,\frac{z}{\epsilon})|\leq \epsilon S$ (see Assumption \ref{regul}.5). Finally, by using the relation $\nu\circ{\bfm \gamma}_\omega^{-1} ={\bfm c}(\omega,z)e^{2{\bfm V}(\omega)}\nu(dz)$, we have
$${\bfm g}_\epsilon(\omega)= \lim_{\beta\to 0}\int_{\{\beta\leq |z|\}}  h(z){\bfm c}(\omega,\frac{z}{\epsilon})e^{2{\bfm V}(\omega)}\nu(dz).$$

We now prove that ${\bfm g}_\epsilon\to 0$ in $L^1(\Omega)$. For each fixed $\beta>0$, Assumption \ref{limitc} implies that the following convergence holds in $L^1(\Omega)$
$$\lim_{\epsilon\to 0}\int_{|z|>\beta}h(z){\bfm c}(\omega,\frac{z}{\epsilon})\nu(dz)={\bfm \theta}(\omega)\int_{|z|>\beta}h(z)\nu(dz)=0.$$ From Lemma \ref{limitz}, the family $\big(\int_{|z|>\beta}h(z){\bfm c}(\omega,\frac{z}{\epsilon})\nu(dz)\big)_{\beta>0}$ converges in $L^\infty(\Omega)$ as $\beta\to 0$ uniformly with respect to $\epsilon>0$. We deduce that, in $L^1(\Omega)$,
$$\lim_{\epsilon\to 0}\lim_{\beta\to 0}\int_{|z|>\beta}h(z){\bfm c}(\omega,\frac{z}{\epsilon})\nu(dz) =0.$$
We conclude the proof with the help of the following estimate
$$\lim_{\epsilon\to 0}\M_\pi\E\Big[\sup_{0\leq t \leq T}\big|\int_0^t{\bfm g}_\epsilon(\tau_{\overline{X}^\epsilon_{r-}})\,dr\big|\Big]\leq \lim_{\epsilon\to 0}\M_\pi\E\Big[\int_0^T|{\bfm g}_\epsilon(\tau_{\overline{X}^\epsilon_{r-}})|\,dr\Big]\leq T\M_\pi[|{\bfm g}_\epsilon|]\to 0$$ as $\epsilon\to 0$.
\qed


\end{document}